\documentclass[preprint,sort&compress,final,12pt]{elsarticle}

\usepackage{amsmath}
\usepackage{amssymb}
\usepackage{graphicx}
\usepackage{latexsym}
\usepackage{epstopdf}
\usepackage{natbib}
\usepackage{color}
\usepackage{hyperref}

\newtheorem{theorem}{Theorem}[section]

\newtheorem{definition}[theorem]{Definition}
\newtheorem{lemma}[theorem]{Lemma}

\numberwithin{equation}{section}

\def\Proof{\noindent{\bf Proof.}~}
\def\qed{\hfill$\square$\smallskip}

\def\dfrac#1#2{\frac{\displaystyle {#1}}{\displaystyle {#2}}}

\journal{\empty}
\date{}

\begin{document}

\begin{frontmatter}

\title{The stability of equilibrium solutions of periodic Hamiltonian systems in the case of degeneracy}

\author[au1,au2]{Nina Xue}

\address[au1]{School of Mathematical Sciences, Beijing Normal University, Beijing 100875, P.R. China.}

\author[au1]{Xiong Li\footnote{ Partially supported by the NSFC (11571041) and the Fundamental Research Funds for the Central Universities. Corresponding author.}}

\address[au2]{School of Mathematics and Information Sciences, Weifang University, Weifang, Shandong, 261061, P.R. China.}

\ead[au2]{xli@bnu.edu.cn}

\begin{abstract}
In this paper we are concerned with the stability of equilibrium solutions of periodic Hamiltonian systems with one degree of freedom in the case of degeneracy, which means that the characteristic exponents of the linearized system have zero real part, and the high order terms must be considered to solve the stability problem. For almost all degenerate cases, sufficient conditions for the stability and instability are obtained.
\end{abstract}


\end{frontmatter}

\section{Introduction}
Consider a Hamiltonian system with one degree of freedom
\begin{equation} \label{intro1}
\dot{x}=\frac{\partial H}{\partial y}(x,y,t),\ \ \ \ \dot{y}=-\frac{\partial H}{\partial x}(x,y,t),
\end{equation}
where $H(0,0,t)=H_x(0,0,t)=H_y(0,0,t)=0$, the dot indicates differentiation with respect to the time $t$. The Hamiltonian function $H$ is continuous and $2\pi$-periodic in $t$, and real analytic in a neighborhood of the origin $(x,y)=(0,0)$, the Taylor series of $H$ in a neighborhood of the origin is assumed to be
\begin{equation} \label{intro2}
H(x,y,t)=H_{2}(x,y,t)+H_{3}(x,y,t)+\cdots +H_{j}(x,y,t)+ \cdots,
\end{equation}
where
\begin{equation} \label{intro3}
H_{j}(x,y,t)=\sum_{\mu+\nu=j}h_{\mu\nu}(t)x^{\mu}y^{\nu},\ \ j=2,3,\cdots,
\end{equation}
and the coefficients $h_{\mu\nu}(t)$ are $2\pi$-periodic with respect to the time $t$.

The question about the stability of an equilibrium position of a Hamiltonian system plays an important role in the knowledge of the qualitative behavior of the system. Also it is well known that for the periodic Hamiltonian system \eqref{intro1},
the problem of knowing about the stability of equilibrium solutions in the sense of Lyapunov is still an open problem.
 Only in some particular cases we have some methods to determine the type of stability. For instance, if there exists a characteristic exponent of the linearized system with non zero real part, the equilibrium solution is unstable. Thus we assume that the linearized system is stable and that the characteristic exponents are pure imaginary, say $\pm i\,\omega\, (\omega\neq 0),$ and the nonlinear part is necessary to answer the question of stability. 
 
For the general elliptic case, Arnold \cite{Arnold} and Moser \cite{Moser}, \cite{Siegel} had obtained some results about stability of the equilibrium solutions under certain nondegeneracy conditions.   Since then, there are plenty
of works about the stability of the trivial solution, one can refer to \cite{Moser1}, \cite{Siegel}
for a detailed description. For recent developments, one may consult \cite{Liu1}, \cite{Ortega3}
and their references therein. For time-periodic Lagrangian equations, an
analytical method called the third order approximation method, had been developed
recently by Ortega in a series of papers \cite{Ortega0}, \cite{Ortega1}, \cite{Ortega2}, \cite{Ortega4}. After that,
some researchers had extended the applications of the third order approximation method,
and some stability results for several types of Lagrangian equations had
been established, one can refer to \cite{Lei}, \cite{Ortega3} for the forced pendulum,
and \cite{Torres1}, \cite{Torres2} for the singular equations.  However, these results can not be applied directly in the degenerate case.

Without loss of generality we assume that a linear $2\pi$-periodic symplectic transformation has already been made such that $H_{2}$ is given by
$$H_{2}(x,y)=\frac{\omega}{2}(x^{2}+y^{2}).$$
Introduce symplectic action angle variables by means of the formula
$$x=\sqrt{2r}\cos\varphi, \ y=\sqrt{2r}\sin\varphi,$$
in the variables $r$ and $\varphi$, the Hamiltonian \eqref{intro2} becomes
\begin{equation} \label{intro4}
H(r,\varphi,t)=H_{2}(r)+H_{3}(r,\varphi,t)+\cdots +H_{j}(r,\varphi,t)+ \cdots,
\end{equation}
here we still denote the new Hamiltonian by $H$,
and
\begin{equation} \label{H_j}
H_{2}(r)=\omega r, \
H_{j}(r,\varphi,t)=\sum_{\mu+\nu=j}h_{\mu\nu}(t)(2r)^\frac{j}{2}\cos^{\mu}\varphi \sin^{\nu}\varphi, \ \  j\geq 3.
\end{equation}

Applying the Lie normal form process, there exists a formal, symplectic, $2\pi$-periodic change of variables which transforms the Hamiltonian in \eqref{intro4} into the following form:

\vspace{0.1cm}
(1). if $\omega$ is an irrational number, then $H=H(r)=\omega r+A_{4}r^{2}+\cdots +A_{2n}r^{n}+\cdots$, i.e., $H$ depends only on the action variable $r$;

\vspace{0.1cm}
(2).  if $\omega=\frac{p}{k}$  is a rational number, where $p, k$ are relatively prime positive integers,
that is, there is a resonance of $k$ order, then $$H=H(r,\varphi)=A_{4}r^{2}+\cdots +A_{2l}r^{l}+H_{k}(r,k\varphi)+\cdots +H_{m}(r,k\varphi)+\cdots,$$
where $2l<k$ is a maximal even number, and $m\geq k$. The details can be found in Lemma \ref{Normal form} of this paper.

When considering the problem of stability, we only need to obtain the normal form up to certain orders, thus the above formal change of variables
can be convergent, and therefore if there exists some $A_{2j}\neq 0$ for $j=2,3,\cdots,$ it follows from Moser twist theorem that the origin in \eqref{intro1} is stable in the Lyapunov sense. Especially, in \cite{Xue} we had proved that when $\omega$ is a Diophantine number and all constants $A_{2j}=0 \,(j=2,3,\cdots)$, then there exists an analytic, symplectic, $2\pi$-periodic change of variables which transforms the Hamiltonian in \eqref{intro4} into the linear part $H_2=\omega r.$ Obviously, in this case, the origin in \eqref{intro1} is also stable. However, if $\omega$ is an irrational but not Diophantine number, the problem about the stability of the origin in \eqref{intro1} still remains open for us.

In the resonance case, that is, $\omega$ is a rational number, there also are many results about the problem of stability, see \cite{Bardin}, \cite{Cabral}, \cite{Mansilla}, \cite{Markeyev1}, \cite{Markeyev2} \cite{Santos1}, \cite{Santos2}, \cite{Santos3} and the references therein. In general, the techniques used in the literature to study the problem about stability of equilibrium solutions are as follows. First, the Lie normal form (the details can be found in Section 2 of this paper) up to a certain order, is obtained. After that, authors often applied Chetaev's theorem or Lyapunov's theorem (see \cite{Meyer}) to prove the instability and Moser twist theorem  \cite{Moser}, \cite{Siegel}) to prove the stability, respectively.
Using this approach, Cabral and Meyer \cite{Cabral} formulated a general stability criterion. This result reads as

\vspace{0.2cm}
\noindent {\bf Theorem A}\, {\it
Let $K(r,\varphi, t)=\psi(\varphi)r^{n}+O(r^{n+\frac{1}{2}})$ as $r\to0^+$, where $n=\frac{m}{2}$ with $m\geq 3$ be an integer. Suppose that $K$ is an analytic function of $\sqrt{r}, \varphi, t$,  and $\tau$-periodic in $\varphi$, $T$-periodic in $t$. If $\psi(\varphi)\neq 0$ for all $\varphi\in \mathbb{R}$, then the origin $r=0$ is a stable equilibrium for the Hamiltonian system
$$\dot{r}=\frac{\partial K}{\partial \varphi}(r,\varphi,t), \ \ \ \dot{\varphi}=-\frac{\partial K}{\partial r}(r,\varphi,t)$$
in the sense that given $\varepsilon>0$, there exists $\delta>0$ such that if $r(0)<\delta$, then the solution is defined for all $t\in\mathbb{R}$ and $r(t)<\varepsilon$. On the other hand, if $\psi(\varphi)$ has a simple zero, i.e., there exists some $\varphi^{\ast}$ such that $\psi(\varphi^{\ast})=0, \psi'(\varphi^{\ast})\neq 0$, then the equilibrium $r=0$ is unstable.}
\vspace{0.2cm}

In this paper, we also consider the case that $\omega=\frac{p}{k}$  is a rational number, that is, there is a resonance of $k$ order, and  $A_{4}=\cdots=A_{2l}=0$ in the normal form (otherwise, the origin is stable). In this case, the Hamiltonian \eqref{intro4} can be transformed into the form in Theorem A. Specially, if $k=2n+1$, i.e., there is a resonance of odd order, the Hamiltonian in \eqref{intro4} can be put into the form
$$H=Br^{\frac{2n+1}{2}}\cos(2n+1)\varphi+\widetilde{H}(r,\varphi, t),$$
which is analytic and $\widetilde{H}=O(r^{n+1})$ is periodic in $t$. According to Theorem A, if $B\neq 0$, then the origin is unstable in the Lyapunov sense. If $B=0$ or $k=2n+2$, i.e., there is a resonance of even order, now the Hamiltonian in \eqref{intro4} has the form
$$H=r^{n+1}\left(C+D\cos2(n+1)\varphi\right)+\widetilde{H}(r,\varphi, t),$$
which is analytic and $\widetilde{H}=O(r^{n+\frac{3}{2}})$ is periodic in $t$. According to Theorem A, if $|C|>|D|$, then the origin is stable in the Lyapunov sense; if $|C|<|D|$, it is unstable. However when $|C|=|D|$, all zeros of the coefficient function  $C+D\cos2(n+1)\varphi$ are two multiplicity, thus Theorem A can not be applied to obtain some information about the type of stability of the origin, which is the case that we will study in the present paper.

Furthermore, we will consider the more general case that all zeros of $\psi(\varphi)$ in Theorem A have multiplicity greater than one. Obviously, Theorem A cannot be applied in this case, which is regarded in the literature as the degenerate case. Hence, in the degenerate case, the terms of order higher than $r^{n}$ in the Hamiltonian $K(r,\varphi,t)$ must be taken into account to solve the stability problem.

There are many results for the stability in the case of degeneracy. As far as we know, Markeyev is the first author to obtain some results on the stability in this case, who in \cite{Markeyev1} and \cite{Markeyev2} studied the stability in the case of fourth resonance for both periodic Hamiltonian systems with one degree of freedom  and autonomous Hamiltonian systems with two degrees of freedom, obtained some sufficient conditions for stability and instability by using Moser twist theorem and Lyapunov theorem, respectively. After that,  Mansilla and Vidal \cite{Mansilla} further investigated the case of even order resonances of periodic Hamiltonian systems with one degree of freedom. Recently, Bardin \cite{Bardin} studied the degenerate case of periodic Hamiltonian systems with one degree of freedom by using Lie normal forms, and established general criteria to solve the stability problem. The basic method of these results is that by considering up to terms of a certain order in normal form, the conditions for stability and instability are then obtained.

Our main results in Theorem \ref{even} and Theorem \ref{odd} (see Section 3) are the generalization of the above results, which enable us to solve the problem about the stability for almost all these degenerate cases. Moreover, compared with \cite{Mansilla} and \cite{Markeyev1}, we formulate more general stability criteria for arbitrary order resonances, and compared with \cite{Bardin}, here we establish the stability criteria which is very easy to verify, and do not need to compute the simple main part in \cite{Bardin} which is very difficult to understand.

The paper is organized as follows. In Section 2, we list some basic definitions and results about resonances and  normal forms that will be useful in our approach. In Section 3, we formulate our main results, that is, Theorem \ref{even} and Theorem \ref{odd}, and give some remarks.  In Section 4, we will prove Theorem \ref{even} and Theorem \ref{odd}. Several examples which can apply our main results are given in Section 5 to illustrate that our results generalize the known ones.

\section{Resonances and normal forms }

We assume that the corresponding linearized system of \eqref{intro1} is stable and the characteristic exponents are $\pm i\,\omega$ with $\omega\in \mathbb{R}$. Let
\begin{equation}\label{nf1}
H^{m}(r,\varphi,t)=H_{2}(r)+H_{3}(r,\varphi,t)+\cdots+H_{m}(r,\varphi,t),
\mathbf{}\end{equation}
where $H_{j}\ (j=2,\cdots, m)$ are given in \eqref{H_j}. Since the resonance relations play an important role in questions of the stability, here we give a precise description for the resonance relationships of system \eqref{intro1}.

\begin{definition}
System \eqref{intro1} is said to possess the resonance relation if there exists an integer $k\neq 0$ such that
$k\omega \in \mathbb{Z}$, and the minimal positive integer number $k$ satisfying $k\omega \in \mathbb{Z}$ is called the order of the resonance. On the other hand, if $k\omega \notin \mathbb{Z}$
for all integers $k$ with $1\leq |k|\leq s$,  we say that system \eqref{intro1} does not present any resonance relations up to the order $s$, inclusively.
\end{definition}

Consider the $\mathbb{Z}$-module
$$M_{\omega}=\{k\in\mathbb{Z}: k\omega\in\mathbb{Z}\}$$
associated to the frequency $\omega$. It is clear that $M_{\omega}=\{0\}$ is equivalent to say that $\omega, 1$ are linearly independent over $\mathbb{Q}$, that is, $M_{\omega}=\{0\}$ if and only if system \eqref{intro1} does not possess any resonance relations. In the opposite case, system \eqref{intro1} possesses the resonance relation and an integer number $k\in M_{\omega}\backslash\{0\}$ is called the number of the resonance and its order of the resonance is $k$, here and hereafter we assume that $k\in M_{\omega}$ is the minimal positive integer number. Then we have the following lemma.
\begin{lemma}\label{Normal form}
Assume that $H^{m}$ ($m$ is a fixed integer greater than 2) is in Lie normal form, if $M_{\omega}=\{0\}$, then $H^{m}=H^{m}(r)$ for all $m$, i.e., $H^{m}$ depends only on the variable $r$; if $M_{\omega}=k\mathbb{Z}\neq\{0\},$ then
$$H^{m}=H^{m}(r,\varphi)=H_{4}(r)+\cdots+ H_{2l}(r)+H_{k}(r,k\varphi)+\cdots +H_{m}(r,k\varphi),$$
where $2l$ is a maximal even number which is less than $k$, and $m\geq k$.
\end{lemma}
\Proof  It follows from \cite{Meyer} that $H^{m}$ is in Lie normal form whenever
\begin{equation} \label{nf2}
\{H_{2},H^{m}\}-\frac{\partial H^{m}}{\partial t}=0,
\end{equation}
where $\{,\}$ is the Poisson bracket.
Since $H_{2}=\omega r$, it follows from \eqref{nf2} that
\begin{equation} \label{nf3}
\omega \frac{\partial H^{m}}{\partial \varphi}=\frac{\partial H^{m}}{\partial t}.
\end{equation}
Assume that
$$ H^{m}(r,\varphi,t)=\sum_{j=2}^{m}\sum_{\mu+\nu=j}h_{\mu\nu}(t)e^{i(\mu-\nu)\varphi}r^{\frac{j}{2}},$$
then equation \eqref{nf3} gives us
$$\sum_{j=2}^{m}\sum_{\mu+\nu=j}i\omega(\mu-\nu) h_{\mu\nu}(t)e^{i(\mu-\nu)\varphi}r^{\frac{j}{2}}
=\sum_{j=2}^{m}\sum_{\mu+\nu=j}h'_{\mu\nu}(t)e^{i(\mu-\nu)\varphi}r^{\frac{j}{2}}.$$
Therefore, the periodic function $h_{\mu\nu}(t)$ satisfies the differential equation
$$ h'_{\mu\nu}(t)=i\omega (\mu-\nu)h_{\mu\nu}(t),$$
whose solutions are
$$h_{\mu\nu}(t)=C_{\mu\nu}e^{i\omega (\mu-\nu)t},\ \  C_{\mu\nu}\in \mathbb{C}.$$

Thus we can rewrite $H^{m}$ as
$$H^{m}(r,\varphi,t)=\sum_{j=2}^{m}\sum_{\mu+\nu=j}C_{\mu\nu}e^{i(\mu-\nu)(\varphi+\omega t)}r^{\frac{j}{2}}.$$
Since $H^{m}$ is a real function in its original variables, by symmetry of subindex $\mu+\nu=j$, then $C_{\mu\nu}=\overline{C}_{\nu\mu}$. Moreover, since $H^{m}$ is of period $2\pi$ in $t$, $H^{m}\neq 0$ if and only if $(\mu-\nu)\omega\in \mathbb{Z},$ i.e., $\mu-\nu\in M_{\omega}$. Thus $\mu=\nu$  if $M_{\omega}=\{0\}$, that is, $H^{m}=H^{m}(r)$; $\mu-\nu\in k\mathbb{Z}$ if $M_{\omega}=k\mathbb{Z}\neq \{0\}$, that is, $H^{m}=H^{m}(r,k(\varphi+\omega t)).$

In the case $M_{\omega}=k\mathbb{Z}\neq \{0\}$, we will make the canonical change of variables
$$ \widetilde{r}=r,  \ \ \widetilde{\varphi}=\varphi+\omega t.$$
Indeed, this transformation can be defined by
$$ r=\frac{\partial S}{\partial \varphi},  \ \ \widetilde{\varphi}=\frac{\partial S}{\partial\widetilde{r}},$$
where the generating function $S$ is
$$S(\widetilde{r},\varphi, t)=\widetilde{r}\cdot(\varphi+\omega t).$$
Since $\frac{\partial S}{\partial t}=\widetilde{r}\cdot \omega,$  the new Hamiltonian has the form
$$\widetilde{H}(\widetilde{r},\widetilde{\varphi})=\widetilde{H}_{4}(\widetilde{r})+\cdots+\widetilde{H}_{2l}(\widetilde{r})
+\widetilde{H}_{k}(\widetilde{r},k\widetilde{\varphi})+\cdots+\widetilde{H}_{m}(\widetilde{r},k\widetilde{\varphi}).$$
Let $\widetilde{H}^{m}=H^{m}, \widetilde{r}=r$ and $\widetilde{\varphi}=\varphi,$  thus the proof of Lemma \ref{Normal form} is completed.\qed

\section{Main results}
We first formulate two general stability criteria for the following Hamiltonian
\begin{equation} \label{intro5}
H(r,\varphi,t)=r^{\alpha}\psi_{0}(\varphi)+r^{\alpha+\gamma}\psi_{1}(\varphi)+O(r^{\alpha+\gamma+\frac{1}{2}}),
\end{equation}
where $\alpha=\frac{m}{2}, \ m\geq 3; \gamma=\frac{n}{2}, \ n\geq 1.$

Our purpose in this paper is to study the stability of the origin of the Hamiltonian \eqref{intro5} in the case of degeneracy, i.e., all zeros of the coefficient function  $\psi_{0}(\varphi)$ have multiplicity greater than one.

Now we are in a position to state our first result. That is, if all zeros of $\psi_{0}(\varphi)$ are even multiplicity, we have the following result.

\begin{theorem}\label{even}
Consider the Hamiltonian system
\begin{equation}\label{th1}
 \dot{\varphi}=\frac{\partial H}{\partial r}(r,\varphi,t), \ \ \dot{r}=-\frac{\partial H}{\partial \varphi}(r,\varphi,t)
 \end{equation}
with the Hamiltonian \eqref{intro5}, and suppose that all zeros of $\psi_{0}(\varphi)$ are even multiplicity. \\
 (A) If every zero $\varphi_{0}$ of $\psi_{0}(\varphi)$ satisfies
 $$\psi_{0}^{(m)}(\varphi_{0})\psi_{1}(\varphi_{0})>0,$$
 where $m>1$ is the multiplicity of the zero $\varphi_{0}.$
 Then the equilibrium $r=0$ of the Hamiltonian system \eqref{th1} is stable in the sense of Lyapunov; \\
 (B) if $\psi_{1}(\varphi)\neq 0$ for all $\varphi$ satisfying $\psi_{0}(\varphi)=0,$ moreover, there exists a zero
$\varphi_{0}$ of $\psi_{0}(\varphi)$ such that
 $$\psi_{0}^{(m)}(\varphi_{0})\psi_{1}(\varphi_{0})<0,$$ where $m>1$ is the multiplicity of the zero $\varphi_{0}.$
 Then the equilibrium $r=0$ of the Hamiltonian system \eqref{th1} is unstable in the sense of Lyapunov.
\end{theorem}

We give two remarks about this result. First of all, if all zeros of $\psi_{0}(\varphi)$ are even multiplicity, there are examples which show that terms whose powers are higher than $\alpha$ can be choose in such a way as to obtain stability or instability, as desired.

For example, consider a system with a Hamiltonian of the form
$$H=(1+\sin\varphi)r^{2}+ar^{3},$$
here we have transformed the Hamiltonian into the normal form, and $a$ is constant coefficient. Obviously, all zeros of $\psi_{0}(\varphi)=1+\sin\varphi$
 have multiplicity two.

 If $a=1$, the equilibrium solution is stable, which can be showed using Lyapunov's theorem on stability, by taking the function $V=H$ as
Lyapunov's function.

If $a=-1$, the equilibrium solution is unstable. In fact, consider the motions at the level $H=0$. At this level either $r=0$ or $r=1+\sin\varphi$.
The first case is of no interest since it corresponds to the equilibrium solution itself. In the second case, i.e., $r=1+\sin\varphi$, the following equations
$$\dot{\varphi}=-(1+\sin\varphi)^{2}, \ \ \dot{r}=-r^{2}\cos\varphi$$
hold. If we put $\varphi(0)=-\frac{\pi}{2}-\mu,$ where $0<\mu\ll 1$, then $r(0)=2\sin^{2}\frac{\mu}{2}\sim \mu^{2}.$ The angle $\varphi$ decreases monotonically with time, as long as it remains in the third quadrant $(-\pi<\varphi<-\frac{\pi}{2})$, and the value of $r$ increases monotonically from small value $r(0)$ to $r=1$, which indicates that the equilibrium solution is unstable.

This example shows that the terms of order higher than $r^\alpha$ in the Hamiltonian \eqref{intro5}, that is, the term $r^{\alpha+\gamma}\psi_{1}(\varphi)$,  must be taken into account to solve the stability problem.

Secondly, if the functions $\psi_0(\varphi)$ and $\psi_1(\varphi)$ have some common zeros, Theorem \ref{even} cannot be applied directly to solve the stability problem of $r=0$. Thus, in this case we have to add the terms of order higher than $r^{\alpha+\gamma}$ into the Hamiltonian \eqref{intro5} to solve the stability problem of $r=0$, which is still a degenerate case, and can use the same idea to study the stability problem of $r=0$.

On the other hand, if $\psi_{0}(\varphi)$ has a zero of odd multiplicity greater than one, our result about the type of stability of the origin is the following.
\begin{theorem}\label{odd}
  If there exists some $\varphi_{0}$ such that
 $$\psi_{0}(\varphi_{0})=\cdots=\psi_{0}^{(m-1)}(\varphi_{0})=0, \ \psi_{0}^{(m)}(\varphi_{0})\neq 0,$$
 and
 \begin{equation*}\label{intro7}
 \psi_{0}^{(m)}(\varphi_{0})\psi'_{1}(\varphi_{0})>0,
 \end{equation*}
 where $m>1$ is odd, then the equilibrium $r=0$ of the Hamiltonian system \eqref{th1} is unstable in the sense of Lyapunov.
\end{theorem}

The proofs of our main results in Theorem \ref{even} and Theorem \ref{odd} will be postponed in Section 4.

Now we return to the Hamiltonian system \eqref{intro1}. Suppose that the characteristic exponents $\pm i\, \omega$ of the linearized system of \eqref{intro1} satisfy a resonance relation, and the order of the resonance is $k$. After introducing action angle variables and the Lie process of normalization up to order $r^{\frac{m}{2}}\, (m>k)$, by Lemma \ref{Normal form}, one can rewrite the Hamiltonian function of system \eqref{intro1} in the following form
\begin{equation}\label{result1}
H=A_4r^2+\cdots+A_{2l}r^l+H_{k}(r,\varphi)+\cdots+H_{m}(r,\varphi)+O(r^{\frac{m+1}{2}}),
\mathbf{}\end{equation}
where $2l$ is a maximal even number which is less than $k$, and $m\geq k$.
Moreover, by Lemma \ref{Normal form} again, $H_s(r,\varphi)\, (k\leq s\leq m)$ in \eqref{result1} has the form
$$H_s(r,\varphi)=\left(A_s+\sum_{j=1}^{[\frac{s}{k}]}(B_{sj}\cos(jk\varphi)+C_{sj}\sin(jk\varphi))\right)r^{\frac{s}{2}}, \ k\leq s\leq m,$$
where $A_s, B_{sj}, C_{sj}$ are some real numbers.

As discussed in Section 1, we may assume that $A_4=\cdots=A_{2l}=0$ in \eqref{result1} (otherwise, the origin is stable). First we consider the case  that $H_k\not\equiv 0$. If $k$ is odd, then $A_k=0$, which is the case that we can apply Theorem A directly,  we may assume that $k$ is even.

Thus, if $k$ is even and $H_k\not\equiv 0$, then
$$H_k(r,\varphi)=(A_k+B_{k1}\cos(k\varphi)+C_{k1}\sin(k\varphi))r^{\frac{k}{2}}.$$
In the degenerate case, that is, $|A_k|=\sqrt{B_{k1}^2+C_{k1}^2}$, we have
$$H_k(r,\varphi)=2A_k\cos^2\frac{k\varphi+\tilde{\varphi}}{2}r^{\frac{k}{2}}:=\psi_0(\varphi)r^{\frac{k}{2}}.$$
It is clear that all zeros of $\psi_{0}(\varphi)=2A_k\cos^2\frac{k\varphi+\tilde{\varphi}}{2}$ are two multiplicity, Theorem \ref{even} can be applied to solve the stability problem of the origin by taking into account the terms of order higher than $r^{\frac{k}{2}}$. Thus, in this case, the stability problem of the origin of system \eqref{intro1} is completely solved by Theorem \ref{even}.

Now we consider the case that $H_k\equiv 0$, which is immediately changed into the discussion for the Hamiltonian \eqref{intro5}. Next we give two particular examples which show that $m>2$ is even in Theorem \ref{even} or $m>1$ is odd in Theorem \ref{odd} can take place.

Firstly, we analyze the case
$$H_k\equiv 0, H_{k+1}\equiv 0, \cdots, H_{3k-1}\equiv 0, H_{3k}\not\equiv0,$$
where $k$ may be even or odd, and $3k\leq m$. In this case, we have
\begin{eqnarray*}
H_{3k}(r,\varphi)&=&\Big(A_{3k}+B_{3k1}\cos k\varphi+C_{3k1}\sin k\varphi+B_{3k2}\cos 2k\varphi+C_{3k2}\sin 2k\varphi \\[0.2cm]
&&+B_{3k3}\cos 3k\varphi+C_{3k3}\cos 3k\varphi\Big)r^{\frac{3k}{2}}.
\end{eqnarray*}
In particular, if
$$A_{3k}=B_{3k2}=C_{3k1}=C_{3k2}=C_{3k3}=0, \ B_{3k1}-3B_{3k3}=0,\ B_{3k3}\not=0,$$
then 
$$H_{3k}(r,\varphi)=4B_{3k3}\cos^3(k\varphi) r^{\frac{3k}{2}}:=\psi_0(\varphi)r^{\frac{3k}{2}}.$$
It is easy to see that all zeros of
$$\psi_0(\varphi)=4B_{3k3}\cos^3k\varphi$$
are three multiplicity. Therefore, the case that $m>1$ is odd in Theorem \ref{odd} takes place.

Similarly we assume 
$$H_k\equiv 0, H_{k+1}\equiv 0, \cdots, H_{4k-1}\equiv 0, H_{4k}\not\equiv0,$$
where $4k\leq m$. Therefore,
\begin{eqnarray*}
H_{4k}(r,\varphi)&=&\Big(A_{4k}+B_{4k1}\cos k\varphi+C_{4k1}\sin k\varphi
+B_{4k2}\cos 2k\varphi+C_{4k2}\sin 2k\varphi \\[0.2cm]
&&+B_{4k3}\cos 3k\varphi+C_{4k3}\cos 3k\varphi
+B_{4k4}\cos 4k\varphi+C_{4k4}\cos 4k\varphi\Big) r^{2k}.
\end{eqnarray*}
In particular, if
$$B_{4k1}=B_{4k3}=C_{4k1}=C_{4k2}=C_{4k3}=C_{4k4}=0, $$
$$ \ A_{4k}-3B_{4k4}=0, \ B_{4k2}-4B_{4k4}=0,\ B_{4k4}\not=0,$$
then
$$H_{4k}(r,\varphi)=8B_{4k4}\cos^4(k\varphi) r^{2k}:=\psi_0(\varphi)r^{2k}.$$
It is easy to see that all zeros of
$$\psi_0(\varphi)=8B_{4k4}\cos^4k\varphi r^{2k}$$
are four multiplicity. Therefore, the case that $m>2$ is even in Theorem \ref{even} takes place.



\section{The proofs of Theorem \ref{even} and Theorem \ref{odd}}
Firstly, we will prove Theorem \ref{even}. The main ideas are the use of Moser twist theorem and Chetaev's theorem to prove the stability and instability of the origin, respectively.

For the sake of convenience, we formulate Chetaev's theorem which will be applied in Theorem \ref{even} to prove the instability of the origin.

\begin{lemma}\label{Chetaev's theorem}
(See[\cite{Meyer}, Chapter 13]) Let $V: O\rightarrow \mathbb{R}$ be a smooth function and $\Omega$ an open subset of $O$ with the following properties:

\noindent (i)$\  \zeta_0\in\partial\Omega$, where $\zeta_0$ is an equilibrium solution of the system $\dot{z}=f(z)$;

\noindent (ii) $ \ V(z)>0$ for $z\in \Omega$;

\noindent (iii) $ \ V(z)=0$ for $z\in \partial\Omega$;

\noindent (iv)$\ \dot{V}(z)=V(z)\cdot f(z)>0$ for $z\in \Omega.$

\noindent Then the equilibrium solution $\zeta_0$ of the system $\dot{z}=f(z)$ is unstable.
\end{lemma}

Now we proceed to prove Theorem \ref{even}.

\vspace{0.2cm}
\noindent{\bf Proof of Theorem \ref{even}}\  In order to prove the item (A), we will use the idea that is similar to \cite{Bardin}, \cite{Mansilla}. Initially we consider the truncated system with the Hamiltonian $H_{0}=r^\alpha\psi_{0}(\varphi)+r^{\alpha+\gamma}\psi_{1}(\varphi).$ By the assumptions of (A) in Theorem \ref{even}, without loss of generality, we may assume that
$$\psi_{0}^{(m)}(\varphi_{0})>0, \ \  \psi_{1}(\varphi_{0})>0,$$
where $m>1$ is the multiplicity of the zero $\varphi_{0}$ of the function $\psi_{0}(\varphi).$
Expanding $\psi_{0}(\varphi)$, $\psi_{1}(\varphi)$ at $\varphi_0$, we obtain
\begin{eqnarray*}\label{pf+1}
H_0&=&r^\alpha\psi_{0}(\varphi)+r^{\alpha+\gamma}\psi_{1}(\varphi) \nonumber \\
&=&r^\alpha\frac{\psi_0^{(m)}(\varphi_0)}{m!}(\varphi-\varphi_0)^m+r^{\alpha+\gamma}\psi_1(\varphi_0)\\
&&+r^\alpha O(\varphi-\varphi_0)^{m+1}+r^{\alpha+\gamma} O(\varphi-\varphi_0), \ |\varphi-\varphi_0|\ll1.
\end{eqnarray*}
Therefore, in this case, $H_{0}$ is positive definite for $r>0$ small enough. Moreover,
$$\frac{d\varphi}{dt}=\frac{\partial H_{0}}{\partial r}=\alpha r^{\alpha-1}\left(\psi_{0}({\varphi})
+\frac{\alpha+\gamma}{\alpha}r^{\gamma}\psi_{1}(\varphi)\right)>0$$
for $r>0$ small enough.

Now we consider the energy equation $H_{0}(r,\varphi)=h_{0}$ for some level of energy $h_{0}$. By the implicit function theorem, we can get $r=r_{0}(\varphi,h_{0}).$ Let $h=h_{0}+\mu, \ |\mu|\ll 1,$  the solution $r(\varphi,h)$ of the equation $H_{0}(r,\varphi)=h$ can be represented by a series of powers of $\mu$,
namely,
$$r(\varphi,h)=r_{0}(\varphi,h_{0})+\mu r_{1}(\varphi,h_{0})+O(\mu^{2})$$
with
$$r_{1}(\varphi,h_{0})=\left(\frac{\partial H_{0}}{\partial r}\right)^{-1}>0.$$

The next step is to introduce the action-angle variables associated to the truncated Hamiltonian $H_{0}$. Denote by $I(h)$ the action variable with $I(0)=0$. It is clear that
$$I=\frac{1}{2\pi}\int^{2\pi}_{0}r(\varphi,h)d\varphi=I_{0}+\mu I_{1}+O(\mu^{2}),$$
where $$I_{0}=\frac{1}{2\pi}\int^{2\pi}_{0}r_{0}(\varphi,h_{0})d\varphi, \ I_{1}=\frac{1}{2\pi}\int^{2\pi}_{0}r_{1}(\varphi,h_{0})d\varphi>0.$$
Since $I_{1}\neq 0$, it follows from the inverse function theorem that the inverse function $h=h(I)$ of $I=I(h)$ exists,  which is an analytic function in a region $0<I\ll 1.$

In the action-angle variables $(I, \theta)$, the Hamiltonian of the full system reads
$$H=h(I)+h_{1}(\theta,I,t),$$
where the function $h_{1}(\theta,I,t)=o(h(I))$ is analytic with respect to $I$ and $\theta$, $2\pi$-periodic in $\theta$ and $T$-periodic in $t$.

According to  Moser twist theorem, the following non-degeneracy condition
\begin{equation}\label{con}
\frac{d^{2}h}{dI^{2}}\neq 0
\end{equation}
guarantees the stability of $I=0$.

Calculations show that
\begin{equation*}\label{con2}
\frac{d^{2}h}{dI^{2}}=\frac{\lambda^{3}}{2\pi}\int^{2\pi}_{0}\left(\frac{\partial H_{0}}{\partial r}\right)^{-3}\frac{\partial^{2}H_{0}}{\partial r^{2}}d\varphi,
\end{equation*}
where $$\lambda=\frac{2\pi}{\int^{2\pi}_{0}\left(\frac{\partial H_{0}}{\partial r}\right)^{-1}d\varphi}.$$
Since $\frac{\partial H_{0}}{\partial r}>0,$ we have $\lambda>0$ for $r>0$ small enough.
By means of direct calculations we have
$$\frac{\partial^{2}H_{0}}{\partial r^{2}}=\alpha(\alpha-1)r^{\alpha-2}\left(\psi_{0}(\varphi)+
\frac{(\alpha+\gamma)(\alpha+\gamma-1)}{\alpha(\alpha-1)}r^{\gamma}\psi_{1}(\varphi)\right)>0$$
for $r>0$ small enough. Thus, the condition \eqref{con} is fulfilled. This completes the proof of the item (A).

Now we prove the item (B). Let
\begin{equation*}
\psi(r,\varphi)=\psi_{0}(\varphi)+r^{\gamma}\psi_{1}(\varphi),
\end{equation*}
then the Hamiltonian \eqref{intro5} reads
$$H(r,\varphi,t)=r^{\alpha}\psi(r,\varphi)+O(r^{\alpha+\gamma+\frac{1}{2}}),$$
where $\alpha=\frac{m}{2}, \ m\geq 3; \gamma=\frac{n}{2}, \ n\geq 1.$

From the hypotheses of the item (B), it follows that for any fixed sufficiently small $r>0$, all real roots of
$\psi(r,\varphi)=\psi_{0}(\varphi)+r^{\gamma}\psi_{1}(\varphi)=0$ are simple. Indeed, choosing any zero $\varphi_0$ of  $\psi_{0}(\varphi)$ with the multiplicity $m$, where $m$ is an even number, expanding $\psi_{0}(\varphi)$, $\psi_{1}(\varphi)$ at $\varphi_0$, we obtain
\begin{eqnarray*}\label{pf+1}
\psi(r,\varphi)&=&\psi_{0}(\varphi)+r^{\gamma}\psi_{1}(\varphi) \nonumber \\
&=&\frac{\psi_0^{(m)}(\varphi_0)}{m!}(\varphi-\varphi_0)^m+O(\varphi-\varphi_0)^{m+1}+r^\gamma\psi_1(\varphi_0)+r^\gamma O(\varphi-\varphi_0).
\end{eqnarray*}

By the assumptions of the item (B), we know that $\psi_1(\varphi_0)\not=0$. If
$$\psi_{0}^{(m)}(\varphi_{0})\psi_{1}(\varphi_{0})>0,$$
then
$$
\psi(r,\varphi)>0 \ (\mbox{or}<0), \ \ \ \ |\varphi-\varphi_0|\ll1,
$$
which implies that $\psi(r,\varphi)=\psi_{0}(\varphi)+r^{\gamma}\psi_{1}(\varphi)=0$ has no any real roots in a sufficiently small neighborhood of $\varphi_{0}$.

If
$$\psi_{0}^{(m)}(\varphi_{0})\psi_{1}(\varphi_{0})<0,$$
then
$$\varphi_{1,2}^*=\varphi_0\pm \left(\frac{-m!r^\gamma \psi_1(\varphi_0)}{\psi_0^{(m)}(\varphi_0)}+O(\varphi-\varphi_0)^{m+1}+r^\gamma O(\varphi-\varphi_0)\right)^{\frac{1}{m}}$$
are two real roots of $\psi(r,\varphi)=0$. Since $\varphi_{1,2}^*\not=\varphi_0$,  and
$$
\psi_{0}(\varphi)=\frac{\psi_0^{(m)}(\varphi_0)}{m!}(\varphi-\varphi_0)^m+O(\varphi-\varphi_0)^{m+1},
$$
then
$$
\psi_{0}'(\varphi_{1,2}^*)\not=0,
$$
which implies that for any fixed sufficiently small $r>0$,
$$
\dfrac{\partial \psi}{\partial \varphi}(r,\varphi_{1,2}^*)\not=0.
$$
Therefore, two real roots $\varphi_{1}^{\ast}$ and $\varphi_{2}^{\ast}$ of $\psi(r,\varphi)=0$ are simple for any fixed sufficiently small $r>0$.

We will use Chetaev's theorem to prove the instability in this case. To this purpose, it is necessary to define a convenient Chetaev's function $V$ and a region $\Omega$ such that  Chetaev's theorem  can be applied, see Lemma \ref{Chetaev's theorem}.

Now we define the Chetaev's function
\begin{equation*}
V=\delta r^{2\alpha+\gamma}-r^{2\alpha}\psi^{2}(r,\varphi),
\end{equation*}
where $\delta$ is a positive number which will be chosen next such that in the region $\Omega$ which we will define soon, we would have that
 $V>0$ in $\Omega$, $V=0$ in $\partial\Omega$ and $\dot{V}>0$ in $\Omega$.

Let
$$
\begin{array}{lll}
D_{a}&=&\{(r,\varphi):V\geq 0,r<a\},\\[0.3cm]
D_{a}^{+}&=&\{(r,\varphi):V\geq 0, r<a, \frac{\partial\psi}{\partial\varphi}(r,\varphi)>0\},\\[0.3cm]
D_{a}^{-}&=&\{(r,\varphi):V\geq 0,r<a, \frac{\partial\psi}{\partial\varphi}(r,\varphi)<0\},
\end{array}
$$
where $a$ is a convenient positive real number. By the hypotheses of the item (B), namely, there exists a zero
$\varphi_{0}$ of $\psi_{0}(\varphi)$ such that
 $$\psi_{0}^{(m)}(\varphi_{0})\psi_{1}(\varphi_{0})<0,$$
 where $m>1$ is the multiplicity of the zero $\varphi_{0},$ we have that, for any fixed sufficiently small $r>0$,
   the equation $\psi(r,\varphi)=0$ has  two simple real roots $\varphi_{1}^{\ast}$ and $\varphi_{2}^{\ast}$ near $\varphi_0$, and we may assume that
   \begin{equation}\label{pf+2}
   \frac{\partial\psi}{\partial\varphi}(r,\varphi_{1}^{\ast})>0, \ \ \ \ \frac{\partial\psi}{\partial\varphi}(r,\varphi_{2}^{\ast})<0.
   \end{equation}
Thus, it is easy to see that $D_{a}\neq\emptyset$,  because
$$V(r,\varphi_{1}^{\ast})=\delta r^{2\alpha+\gamma}-r^{2\alpha}\psi^{2}(r,\varphi_{1}^{\ast})
=\delta r^{2\alpha+\gamma}\geq 0.$$
 On the other hand, from \eqref{pf+2}, it follows that $D_{a}^{+}$ and $D_{a}^{-}$ are nonempty subset of $D_{a}$.

In $D_{a}$, it holds that
\begin{equation*}
\delta r^{2\alpha+\gamma}-r^{2\alpha}\psi^{2}(r,\varphi)\geq 0,
\end{equation*}
i.e.,

\begin{equation}\label{result2}
\psi^{2}(r,\varphi)\leq \delta r^{\gamma}\leq \delta a^{\gamma}<\delta
\end{equation}
for $0<a<1.$
Now we will obtain the conditions on $\delta$, such that, for $r>0$ sufficiently small,
$ D_{a}=D_{a}^{+}\cup D_{a}^{-}$  and $\overline{D_{a}^{+}}\cup \overline{D_{a}^{-}}=\{r=0\}$, where the overline on the sets stands for their closure. To do this, it is sufficient to show that
$$\frac{\partial\psi}{\partial\varphi}(r,\varphi)\neq 0$$
in $D_{a}$ for every $r>0$ sufficiently small. In fact,
let us assume the opposite and we will get a contradiction.
We will assume that for $r=\widetilde{r}>0$ sufficiently small, there exists some $\varphi=\widetilde{\varphi}$ such that $$\frac{\partial\psi}{\partial\varphi}(\widetilde{r},\widetilde{\varphi})= 0$$
Since all real roots of the equation $\psi(r,\varphi)=0$ are simple, we have that $\psi(\widetilde{r},\widetilde{\varphi})\neq 0. $
Take $0< \delta\leq\psi^2(\widetilde{r},\widetilde{\varphi})$, thus, by \eqref{result2} we get
$$\delta\leq\psi^2(\widetilde{r},\widetilde{\varphi})<\delta,$$
which is a contradiction. Therefore, $D_{a}=D_{a}^{+}\cup D_{a}^{-}.$

If for any $r>0$ and $\varphi \in \mathbb{R}$, we have $\frac{\partial\psi}{\partial\varphi}(r,\varphi)\neq 0.$  Then we can choose an arbitrary positive number $\delta$. Thus, in this case, we also have
$D_{a}=D_{a}^{+}\cup D_{a}^{-}.$  On the other hand, it is clear that $\overline{D_{a}^{+}}\cup \overline{D_{a}^{-}}=\{r=0\}$.

Our next step is to calculate the derivative of $V$ through the solutions of the Hamiltonian system associated to $H$.
Calculations show that
\begin{eqnarray*}
\dot{V}&=&\{V,H\}=\{\delta r^{2\alpha+\gamma}-r^{2\alpha}\psi^{2}(r,\varphi),H\} \\[0.4cm]
&=&\delta(2\alpha+\gamma)r^{2\alpha+\gamma-1}\{r,H\}
-2r^{\alpha}\psi(r,\varphi)\{r^{\alpha}\psi(r,\varphi),H\}\\[0.2cm]
&=&-\delta(2\alpha+\gamma)r^{2\alpha+\gamma-1}\frac{\partial H}{\partial\varphi}
-2r^{\alpha}\psi(r,\varphi)\{r^{\alpha}\psi(r,\varphi),O(r^{\alpha+\gamma+\frac{1}{2}})\} \\[0.2cm]
&=&-\delta(2\alpha+\gamma)r^{2\alpha+\gamma-1}\left(r^{\alpha}\frac{\partial\psi}{\partial\varphi}(r,\varphi)
+O(r^{\alpha+\gamma+\frac{1}{2}})\right)\\[0.2cm]
&&{}-2r^{\alpha}\psi(r,\varphi)O(r^{2\alpha+\gamma-\frac{1}{2}})\\[0.2cm]
&=&- \delta(2\alpha+\gamma)r^{3\alpha+\gamma-1}\frac{\partial\psi}{\partial\varphi}(r,\varphi)+
O(r^{3\alpha+\gamma-\frac{1}{2}}),
\end{eqnarray*}
where $\dot{V}$ means the derivative of $V$ through the solutions.
Obviously, $\dot{V}>0$ in $D_{a}^{-}$. For application of Chetaev's theorem we consider only one connected component $\Omega$ in $D_{a}^{-}$.
 Thus $V$ is a Chetaev's function in the region $\Omega$, the Chetaev's theorem guarantees the instability of $r=0$.

 Thus, Theorem \ref{even} is proved completely.
\qed

Now we give the proof of Theorem \ref{odd}.

\vspace{0.2cm}
\noindent {\bf Proof of Theorem \ref{odd}}  By the hypothesis of Theorem \ref{odd}, without loss of generality, we assume that
 \begin{equation}\label{+1}
\psi_{0}^{(m)}(\varphi_{0})<0, \ \psi'_{1}(\varphi_{0})<0.
\end{equation}
 Expanding $\psi_{0}(\varphi)$ at $\varphi_{0}$ yields that
\begin{equation}\label{B1}
\psi_{0}(\varphi)=\frac{\psi_{0}^{(m)}(\varphi_{0})}{m!}(\varphi-\varphi_{0})^{m}
+O(\varphi-\varphi_{0})^{m+1},
\end{equation}
where $m>1$ is the multiplicity of the zero $\varphi_{0}$ of the function $\psi_{0}(\varphi)$, and $m$ is odd.

From \eqref{+1} and \eqref{B1}, there exists  a number $0<\delta_{1}<1$ such that
\begin{eqnarray*}
\psi'_{0}(\varphi)=\frac{\psi_{0}^{(m)}(\varphi_{0})}{(m-1)!}(\varphi-\varphi_{0})^{m-1}
+O(\varphi-\varphi_{0})^{m}<0
\end{eqnarray*}
holds for any $\varphi$ satisfying $0<|\varphi-\varphi_{0}|<\delta_{1}$. Since $\psi'_{1}(\varphi_{0})<0,$ there exists $\delta_{2}>0$ such that
$\psi'_{1}(\varphi)<0, $
for any $\varphi$ satisfying $ |\varphi-\varphi_{0}|<\delta_{2}.$

Let $\delta=\min\{\delta_{1}, \delta_{2}\}$, and
$$ S=\left\{(\varphi,r):|\varphi-\varphi_{0}|<\delta, \  r>0 \ \ \mbox{small} \ \ \mbox{enough}\right\}.$$
In $S$, it is easy to see that
\begin{eqnarray}\label{B3}
\dot{r}&=&-H_{\varphi}=-r^{\alpha}\psi'_{0}(\varphi)-r^{\alpha+\gamma}\psi'_{1}(\varphi)
+O(r^{\alpha+\gamma+\frac{1}{2}})\nonumber\\[0.2cm]
&=&-\frac{\psi_{0}^{(m)}(\varphi_{0})}{(m-1)!}(\varphi-\varphi_{0})^{m-1}r^{\alpha}
-r^{\alpha+\gamma}\psi'_{1}(\varphi)\nonumber\\[0.2cm]
&&{}+O(r^{\alpha+\gamma+\frac{1}{2}})+O(\varphi-\varphi_{0})^{m}>0,
\end{eqnarray}
and
\begin{eqnarray}\label{B4}
\dot{\varphi}&=&H_{r}=\alpha r^{\alpha-1}\psi_{0}(\varphi)
+O(r^{\alpha-\frac{1}{2}})\nonumber\\[0.2cm]
&=&\alpha r^{\alpha-1}\left(\psi_{0}(\varphi)-\psi_{0}(\varphi_{0})\right)+O(r^{\alpha-\frac{1}{2}})\nonumber\\[0.2cm]
&=&\alpha r^{\alpha-1}\psi'_{0}(\varphi_{1})(\varphi-\varphi_{0})+O(r^{\alpha-\frac{1}{2}}),
\end{eqnarray}
where $\varphi_1$ satisfies $0<|\varphi_1-\varphi_0|<\delta_1.$

Since $\psi'_{0}(\varphi_{1})<0,$ from \eqref{B3} and \eqref{B4}, we have that any solution $\left(\varphi(t),r(t)\right)$ of the Hamiltonian system \eqref{th1} with the initial condition $(\varphi(0),r(0))\in S$ always stay in $S$ and satisfies $r'(t)>0$ for all $t\geq 0$. Thus $r=0$ is unstable in the sense of Lyapunov.
\qed

\section{Some examples}
Our main results extends and generalizes several results existing in the literature. For example, a critical case of fourth order resonance studied by Markeyev \cite{Markeyev1} can apply our results to obtain the stability and instability. In such case, the Hamiltonian function is reduced to the following
$$H=(1-\cos4\varphi)r^{2}+\left(s+k(1-\cos4\varphi)\right)r^{3}+O(r^{4}),$$
where $s=1$ or $s=-1$ and $k\in\mathbb{R}$. He proved that, if $s=1$, the equilibrium position $r=0$ is stable, if $s=-1$, it is unstable.

However, by Theorem \ref{even}, we note that in this case the functions $\psi_{0}(\varphi)$ and $\psi_{1}(\varphi)$ are given by
$$\psi_{0}(\varphi)=1-\cos4\varphi, \ \ \psi_{1}(\varphi)=s+k(1-\cos4\varphi).$$
Obviously, all zeros of $\psi_{0}(\varphi)$ are $\varphi_{j}=\frac{j\pi}{2}, j\in \mathbb{Z}$, they all have multiplicity two. Thus, if $s=1$, we have
$$\frac{d^{2}\psi_{0}}{d\varphi^{2}}(\varphi_{j})\psi_{1}(\varphi_{j})>0,$$
according to the item (A) of Theorem \ref{even}, $r=0$ is stable.
On the contrary, if $s=-1$, we get
$$\frac{d^{2}\psi_{0}}{d\varphi^{2}}(\varphi_{j})\psi_{1}(\varphi_{j})<0,$$
by the item (B) of Theorem \ref{even}, $r=0$ is unstable.

Another result that we find in the literature concerns with the case of even order resonance for periodic Hamiltonian system with one degree of freedom, which had been proved by Mansilla and Vidal \cite{Mansilla} using the approach of Markeyev \cite{Markeyev1}. In this case,  the Hamiltonian function is reduced to the following
$$H=(1-\cos q\varphi)r^{\frac{q}{2}}+\left(s+k(1-\cos q\varphi)\right)r^{\frac{q}{2}+1}+O(r^{\frac{q+3}{2}}),$$
where the order $q>2$ of resonance is even, $s=1$ or $s=-1$ and $k\in\mathbb{R}$. They proved that, if $s=1$, the equilibrium position $r=0$ is stable, if $s=-1$, it is unstable.

However, by Theorem \ref{even}, we note that in this case the functions $\psi_{0}(\varphi)$ and $\psi_{1}(\varphi)$ are given by
$$\psi_{0}(\varphi)=1-\cos q\varphi, \ \ \psi_{1}(\varphi)=s+k(1-\cos q\varphi).$$
Obviously, all zeros of $\psi_{0}(\varphi)$ are $\varphi_{j}=\frac{2j\pi}{q}, j\in \mathbb{Z}$, they all have multiplicity two. Thus, if $s=1$, we obtain
$$\frac{d^{2}\psi_{0}}{d\varphi^{2}}(\varphi_{j})\psi_{1}(\varphi_{j})>0,$$
according to the item (A) of Theorem \ref{even}, $r=0$ is stable.
On the contrary, if $s=-1$,  we have
$$\frac{d^{2}\psi_{0}}{d\varphi^{2}}(\varphi_{j})\psi_{1}(\varphi_{j})<0,$$
by the item (B) of Theorem \ref{even}, $r=0$ is unstable.


\section*{References}
\bibliographystyle{elsarticle-num}

\end{document}